%% file: distributions.tex
\documentclass{amsart}
\include{canonicalheader}
\renewenvironment{proof}{\paragraph{\textit{Démonstration.}}}{\hfill$\square$}
\newcommand{\Lo}{\Log_{a_p}}
\usepackage{amsthm}
\usepackage{mathtools}  
\usepackage[frenchb, english]{babel}
\usepackage{fancyhdr}
\usepackage{bbm}

\begin{document}

\title{La matrice de logarithme en termes de chiffres $p$-adiques\\
(The logarithm matrix in terms of $p$-adic digits)}
\lhead{Florian Sprung}
\rhead{La matrice \`{a} logarithme en termes de chiffres $p$-adiques}
\selectlanguage{french} 
\begin{abstract}
Nous donnons une nouvelle description de la matrice de logarithme d'une forme modulaire en termes de distributions, généralisant le travail de Dion et Lei pour le cas $a_p=0$. Ce qui nous permet d’inclure le cas $a_p\neq0$ est une nouvelle définition, celle d'une matrice de distributions, et la caractérisation de cette matrice par de chiffres $p$-adiques. On peut appliquer ces méthodes au cas correspondant d'une distribution à plusieurs variables.

(We give a new description of the logarithm matrix of a modular form in terms of distributions, generalizing the work of Dion and Lei for the case $a_p=0$. What allows us to include the case $a_p\neq0$ is a new definition, that of a distribution matrix, and the characterization of this matrix by $p$-adic digits. One can apply these methods to the corresponding case of distributions in multiple variables.)
\end{abstract}

\author{Florian Sprung}
\address{Florian Sprung, 
 School of Mathematical and Statistical Sciences\\
Arizona State University\\
Tempe, AZ 85287-1804\\ USA}
\email{florian.sprung@asu.edu}

\maketitle

\section{Introduction}
Soit $f=\sum_{n\geq 1}^{\infty} a_nq^n$ une forme parabolique de poids $2$ propre pour les opérateurs de Hecke telle que $a_1=1$ et de Nebentypus $\epsilon$. Soit $p$ un nombre premier tel que $p$ ne divise pas le conducteur de $f$. Par exemple, on peut prendre pour $f$ la forme associée à une courbe elliptique telle qu'elle ait bonne réduction en $p$ (ordinaire ou supersingulière). Une construction due à Mazur-- Swinnerton-Dyer, Amice--Vélu et Vi\v{s}ik \cite{mazurmc,amicevelu,vishik} nous permet d'associer des fonctions $L$ $p$-adiques $L_\lambda$ pour chaque racine $\lambda$ du polynôme $X^2-a_pX+\epsilon(p)p$ de Hecke telle que $v_p(\lambda)<1$, où $v_p$ et la valuation $p$-adique telle que $v_p(p)=1$. 

Dans le cas ordinaire, où $v_p(\lambda)=0$, on peut utiliser $L_\lambda$ pour formuler une conjecture principale d'Iwasawa, voir par exemple \cite{mazurmc}. Dans le cas supersingulier, où $v_p(\lambda)>0$, la situation n'est pas idéale: Il y a deux choix $\alpha$ et $\beta$ pour la racine $\lambda$, et les deux fonctions $L_\alpha$ et $L_\beta$ sont analytiques, mais ne sont pas bornées sur le disque compact $\Zp$ comme dans le cas ordinaire. Néanmoins, il existe une factorisation de type distribution$ = $mesure $\times$ distribution
\[(L_\alpha,L_\beta)=(L_\sharp,L_\flat)\Lo\]
où $\Lo$ est une matrice de dimension $2\times 2$. Comme les deux fonctions $L_\sharp$ et $L_\flat$ sont des fonctions analytiques bornées sur $\Z_p$, elles sont plus convenables pour formuler des conjectures principales, voir \cite{koba} pour le cas $a_p=0$ et \cite{shuron} pour le cas général pour les courbes elliptiques.

Les coefficients de $\Lo$ (la matrice de logarithme) sont des distributions (après avoir identifié ces fonctions analytiques avec leurs distributions, voir \cite[Théorème 1.3]{colmez})de mêmes ordres que celles de $L_\alpha$ et $L_\beta$. Dans le cas $a_p=0$, les coefficients de $\Lo$ sont des produits infinis de polynômes cyclotomiques, les logarithmes $p$-adiques signés dus à Pollack \cite{pollack}. Dion et Lei \cite{DL} ont donné une description concrète des distributions dont la transformée d'Amice--Mahler sont ces logarithmes signés. 

Dans le cas général, on peut définir $\Lo$ comme un produit infini de matrices \cite{ant} qui contient les polynômes cyclotomiques. Dans le cas ordinaire, seulement les termes de la première colonne de $\Lo$ convergent.

Le résultat principal généralise le résultat de Dion--Lei au cas général et nous fournit une caractérisation de la matrice $\Lo$ par des chiffres $p$-adiques.

On appelle une matrice $\smat{\mu_{11} & \mu_{12}\\ \mu_{21} & \mu_{22}}$ de distributions $\mu_{11}$, etc. une \textit{matrice de distributions}. Pour simplifier la version du théorème dans l'introduction, on suppose ici que $\epsilon(p)=1$.

\begin{thé} Soit $b\in\Zp$ tel que \[b\equiv b_0+b_1p^1+\cdots b_{n-1}p^{n-1}\pmod{p^n}\text{ avec }b_i\in\{0,\cdots p-1\}.\] Soit $\mu_{a_p}$ la matrice des distributions dont sa transformée est $\Lo$. (Dans le cas ordinaire, on ignore simplement la deuxième colonne.)
	
	On dénote par $m_1,\cdots, m_l$ les longueurs des chaines des chiffres $=0$:
	
	$$(b_0,\cdots b_{n-1})=(\underbrace{0,\cdots, 0}_{m_1},\neq 0, \underbrace{0,\cdots, 0}_{m_2},\neq 0,\cdots ,\neq 0,\underbrace{0,\cdots, 0}_{m_l})$$
	
	Alors la matrice $\mu_{a_p}$ est caractérisée par
	
 \[\mu_{a_p}(b+p^n\Zp)=	\smat{ a_p & 1 \\ -1 & 0}^{m_1}\smat{0 & 0 \\ -1 & 0} \smat{ a_p & 1 \\ -1 & 0}^{m_2}\smat{0 & 0 \\ -1 & 0} \cdots \smat{0 & 0 \\ -1 & 0} \smat{ a_p & 1 \\ -1 & 0}^{m_l}R_n,\] 
c.-à-d. la matrice $\smat{0 & 0 \\ -1 & 0}$ représente un chiffre non nul et la matrice $\smat{ a_p & 1 \\ -1 & 0}$ représente le chiffre $0$. (Dans le cas ordinaire, cela caractérise la première colonne, comme on a ignoré la deuxième.)
\end{thé}

Pour la définition de $R_n$, voir Section \ref{logarithmmatrix}. Si $p$ est impair, on a $R_n=\frac{1}{p^{2+n}}\smat{-\beta^{n+2} & -\alpha^{n+2}\\ \beta^{n+{3+n}} & \alpha^{n+3}}$.

Quelques observations sur le produit des matrices:
 \begin{enumerate}
 	\item On a $\mu_{a_p}(b+p^n\Zp)=\smat{0 & 0 \\ 0 & 0}$ si on a deux chiffres non nuls consécutifs.
 	\item Si $a_p=0$, on a $\smat{0 & 0 \\ -1 & 0} \smat{ a_p & 1 \\ -1 & 0}^{m_i}\smat{0 & 0 \\ -1 & 0} =\smat{0 & 0 \\ 0 & 0}$ pour $m_i$ pair. Par conséquent, si un nombre pair de $0$s sépare deux chiffres non nuls, on a automatiquement $\mu_{a_p}(b+p^n\Zp)=\smat{0 & 0 \\ 0 & 0}$. Mais si $m_2,\cdots,m_l$ sont tous impairs, le produit n'est pas $\smat{0 & 0 \\ 0 & 0}$. (On obtient $\smat{0 & \pm1 \\ 0 & 0}R_n$ ou $\smat{0 & 0 \\ 0 & \pm1}R_n$.) C'est-à-dire, $\mu_{a_p}(b+p^n\Zp)$ peut être non nulle si et seulement si touts les chiffres impairs de $b$ sont $0$ ou tous les chiffres pairs le sont. Cela explique la définition des ensembles $S_n^\pm$ dans \cite{DL}.
\end{enumerate}

Il y existe divers analogues de la matrice de logarithme qui jouent un rôle important dans les travaux de Lei, Loeffler, Zerbes et Büyükboduk (mais la construction est différente -- ils utilisent la théorie de modules de Wach). Voir par exemple \cite{ponsinet} pour une théorie récente. La question d'étudier leurs matrices de distributions provenant de la théorie des modules de Wach semble très intéressante.

\section{D\'{e}finition et propri\'{e}t\'{e}s de la matrice \`{a} logarithme}\label{logarithmmatrix}
Nous rappelons la définition de la matrice de logarithme définie dans \cite{ant}.
Soit $\Phi_{p^n}(X)=\sum_{i\geq0}^{p-1} X^{ip^{n-1}}$ le $p^n$-ième polynôme cyclotomique. 
Nous dénotons $\epsilon(p)=\pm 1$ simplement par $\varepsilon$. On pose

$$\Log_{a_p}^{(n)}(1+T):=\smat{a_p &1\\ - \varepsilon\Phi_p(1+T) & 0}\smat{a_p & 1\\ -\varepsilon\Phi_{p^2}(1+T) & 0}\cdots \smat{a_p & 1\\ -\varepsilon\Phi_{p^n}(1+T) & 0} R_n $$
où $R_n$ dénote la matrice de racines. Cette matrice dépend des racines $\alpha$ et $\beta$ du polynôme de Hecke $X^2-a_pX+\varepsilon p=0$ que l'on choisit telles que $|\alpha|\geq|\beta|$ pour une valeur absolue $p$-adique quelconque. 
La matrice $R_n$ est définie par
 $$R_n:=\smat{a_p & 1\\ -\varepsilon p & 0}^{-2-n}\smat{-1 & -1\\ \beta&\alpha}=\frac{1}{p^{2+n}}\smat{-\beta^{2+n} & -\alpha^{2+n}\\ \beta^{3+n} & \alpha^{3+n}}$$ si $p\geq 3$, et $$R_n:=\smat{a_p & 1\\ -\varepsilon p & 0}^{-3-n}\smat{-1 & -1\\ \beta&\alpha}=\frac{1}{p^{3+n}}\smat{-\beta^{3+n} & -\alpha^{3+n}\\ \beta^{4+n} & \alpha^{4+n}}$$ si $p=2$. \footnote{Il nous semble qu'il y a une très petite erreur dans \cite{DL}:  Dans \cite{DL}, touts les résultats nécessitent l'hypothèse $p>2$. (La première formule \cite[Formule 1]{DL} est incorrecte si $p=2$. La formule correcte de Pollack  est $L_\lambda=\log_p^-L_++\frac{1}{2}\log_p^+L_-$ si $p=2$, voir \cite[Theorem 5.6]{pollack}.)}

\begin{déf}
	La matrice de logarithme est $\Log_{a_p}(1+T):=\lim_{n\rightarrow \infty} \Log_{a_p}^{(n)}(1+T)$.
\end{déf}

\begin{lemme}\label{evaluation}
	Les termes de $\Log_{a_p}(1+T)$ convergent vers des séries dans $\Q_p(\alpha)[[T]]$ si $|a_p|>0$ (cas supersingulier) et on peut dire de même sur les termes dans la première colonne si $|a_p|=0$. De plus, on sait évaluer $\Log_{a_p}(T)$ au point $T=\zeta_{k}-1$, où $\zeta_k$ dénote une racine $p^k$-ième primitive:
	$$\Log_{a_p}(\zeta_k)=\Log_{a_p}^{(k)}(\zeta_k)$$
\end{lemme}
\begin{proof}
	\cite[Lemma 4.8]{ant} ou \cite[Lemma 4.4]{shuron}; l'observation clé est que $\Phi_{p^i}(\zeta_k)=p$ si $i>k$.
\end{proof}

\section{Matrices à distribution}
	\begin{déf}
		La \textbf{transformée} (d'Amice--Mahler) d'une distribution $\mu$ sur $\Z_p$ est 
		\[A_\mu(T)=\int_{\Z_p}(1+T)^x\mu(x).\]
	\end{déf}

	\begin{déf}
		Étant donné $n\times n$ distributions $\mu_{11},\mu_{12},\cdots,\mu_{1n},\mu_{21},\cdots,\mu_{nn}$, la \textbf{matrice de distributions} est la matrice
		\[\mu:=\smat{\mu_{11} & \mu_{12} & \cdots & \mu_{1n}\\
			\mu_{21}& \cdots & &\mu_{2n}\\
		\vdots & &\ddots & \vdots\\
	\mu_{n1} & \cdots &\cdots& \mu_{nn}}.\]
	\end{déf}	

	\begin{déf}
		La \textbf{transformée} d'une matrice de distributions $\mu$ comme ci-dessus est
		\[A_\mu(T):=\int_{\Zp}(1+T)^x\mu(x):=\smat{\int_{\Zp}(1+T)^x\mu_{11}(x) & \int_{\Zp}(1+T)^x\mu_{12}(x) & \cdots &\int_{\Zp}(1+T)^x \mu_{1n}(x)\\
		\int_{\Zp}(1+T)^x	\mu_{21}(x)& \cdots & &\int_{\Zp}(1+T)^x\mu_{2n}(x)\\
		\vdots & &\ddots & \vdots\\
		\int_{\Zp}(1+T)^x\mu_{n1}(x) & \cdots &\cdots& \int_{\Zp}(1+T)^x\mu_{nn}(x)}.\]
	\end{déf}

\section{Le th\'{e}or\`{e}me}

\begin{thé} Soit $b\in\Zp$ tel que \[b\equiv b_0+b_1p^1+\cdots b_{n-1}p^{n-1}\pmod{p^n}\text{ avec }b_i\in\{0,\cdots p-1\}.\] Soit $\mu_{a_p}$ la matrice des distributions dont sa transformée est $\Lo$. (Dans le cas ordinaire, on ignore simplement la deuxième colonne.)
	
	On dénote par $m_1,\cdots, m_l$ les longueurs des chaines des chiffres $=0$:
	
	$$(b_0,\cdots b_{n-1})=(\underbrace{0,\cdots, 0}_{m_1},\neq 0, \underbrace{0,\cdots, 0}_{m_2},\neq 0,\cdots ,\neq 0,\underbrace{0,\cdots, 0}_{m_l})$$
	
	Alors on a
	
	\[\mu_{a_p}(b+p^n\Zp)=	\smat{ a_p & 1 \\ -\epsilon & 0}^{m_1}\smat{0 & 0 \\ -\epsilon & 0} \smat{ a_p & 1 \\ -\epsilon & 0}^{m_2}\smat{0 & 0 \\ -\epsilon & 0} \cdots \smat{0 & 0 \\ -\epsilon & 0} \smat{ a_p & 1 \\ -\epsilon & 0}^{m_l}R_n,\] 
	c.-à-d. la matrice $\smat{0 & 0 \\ -\epsilon & 0}$ représente un chiffre non nul et la matrice $\smat{ a_p & 1 \\ -\epsilon & 0}$ représente le chiffre $0$. (Dans le cas ordinaire, on ignore la deuxième colonne après avoir calculé le produit des matrices.)
	
	De plus, la matrice $\mu_{a_p}$ est caractérisée par ce théorème si $|a_p|>0$, c.-à-d. dans le cas supersingulier. Dans le cas ordinaire, on peut dire de même sur les termes de la première colonne de la matrice $\mu_{a_p}$.
\end{thé}
\begin{proof}
	Soit $\mathbbm{1}_{b+p^n\Zp(x)}$ la fonction caractéristique de $b+p^n\Zp$. On a  $$\mathbbm{1}_{b+p^n\Zp(x)}=\frac{1}{p^n}\sum_{\zeta\in\mu_{p^n}}\zeta^{x-b}.$$ Par conséquent,$$\mu_{a_p}(b+p^n\Zp)=\int_{\Zp}\mathbbm{1}_{b+p^n\Zp}(x)\mu_{a_p}(x)=\sum_{\zeta\in\mu_{p^n}}\frac{1}{p^n}\int_{\Zp}\zeta^{x-b}\mu_{a_p}(x)$$
	
	$$=\frac{1}{p^n}\sum_{\zeta\in\mu_{p^n}}\zeta^{-b}\int_{\Zp}\zeta^x\mu_{a_p}(x)=\frac{1}{p^n}\sum_{\zeta\in\mu_{p^n}}\zeta^{-b}\Lo(\zeta).$$
	Pour calculer cette somme, on utilise la Proposition \ref{Proposition} ci-dessous. La caractérisation suit du fait que les entrées dans $\Lo$ (dans le cas supersingulier) et les entrées dans la première colonne de $\Lo$ (cas ordinaire) sont des distributions d'ordre $o(\log_p(1+T))$, voir \cite[Proposition 4.20]{ant}.
\end{proof}


\begin{proposition}\label{Proposition}  Soit $b=b_0+b_1p^1+\cdots b_{n-1}p^{n-1}$ avec $b_i \in \{0,\cdots, p-1\}$. Pour chaque $b_i$, on définit (sa représentation chromatique) par $$Y_i=\begin{cases} \smat{ a_p & 1 \\ -\varepsilon & 0} & \text{si $b_i=0$ et }\\\smat{0 & 0 \\ 0 & 1} \smat{ a_p & 1 \\ -\varepsilon & 0}=\smat{0 & 0 \\ -\varepsilon & 0} & \text{si $b_i\neq 0$.}\end{cases}$$
	
	Alors 
	
	$$ \sum_{\zeta\in\mu_{p^n}}\zeta^{-b}\Log_{a_p}(\zeta)=p^n Y_0Y_1\cdots Y_{n-1}R_n.$$
	
\end{proposition}

\section{D\'{e}monstration de la proposition}

\begin{lemme}\label{constants}
	Soit \[P(x)=\sum_{i\geq -(p^n-1)}^{i \leq p^n-1}c_ix^i \in \C_p[x,x^{-1}],\]
	où $\C_p$ dénote les nombres complexes $p$-adiques. On a alors
	
	$$\sum_{\zeta \in \mu_{p^n} }P(\zeta)=c_0p^n.$$
\end{lemme}
\begin{proof}
	On a $\sum_{\zeta \in \mu_{p^n}} \zeta^k=\begin{cases} p^n & \text{si $k\equiv 0 \mod{p^n}$}\\ 0 & \text{sinon.}\end{cases}$
\end{proof}
\renewenvironment{proof}{\paragraph{\textit{Démonstration de la Proposition.}}}{\hfill$\square$}
\begin{proof}
	
	Il suit du Lemme \ref{evaluation} que $$\zeta^{-b}\Lo(\zeta)= \zeta^{-b}\Lo^{(n)}(\zeta),$$ alors chaque terme est de la forme $P(\zeta)$ pour des polynômes de Laurent $P(x)$ comme dans le Lemme \ref{constants} (où l'on autorise les puissances négatives). En appliquant ce lemme \ref{constants}, on voit que 
	$$ \sum_{\zeta\in\mu_{p^n}}\zeta^{-b}\Log_{a_p}(\zeta)=p^n\times \text{matrice des termes constants dans $x^{-b}\Lo^{(n)}(x)$}.$$
	
	Montrons alors que les termes constants sont donnés comme dans l'énoncé de la proposition. Posons  $$\phi_i(x):=\smat{1 & 0 \\ 0 & \Phi_{p^i}(x)} \text{ et }  Y:=\smat{ a_p & 1 \\ -\varepsilon & 0}.$$
	Par définition, on a$$\Lo^{(n)}(x)=\phi_0(x)Y\phi_1(x)\cdots Y \phi_{n-1}(x)YR_n,$$
	
	et aussi $$x^{-b}=x^{-b_0}x^{-b_1p}\cdots x^{-b_{n-1}p^{n-1}}.$$

	En combinant ces deux observations, on a 
	
	$$x^{-b}\Lo^{(n)}(x)=x^{-b_0}\phi_0(x)Yx^{-b_1p}\phi_1(x)\cdots Yx^{-b_{n-1}p^{n-1}} \phi_{n-1}(x)YR_n.$$
	
	En notant que $$x^{-b_i}\phi_{i+1}=\smat{ x^{-b_i} & 0 \\ 0 & \sum_{k\geq0}^{p-1}x^{(k-b_i)p^i}},$$
	
	on voit que les termes dans $x^{-b}\Lo^{(n)}(x)$ sont des combinaisons linéaires des termes de la forme $ x^{j_ip^i}$ avec $j_i\in\{-(p-1),\cdots,0,\cdots(p-1)\}$.
	
	Mais on cherche les termes constants, c.-à-d. la contribution des (multiples des) termes $x^{j_ip^i}$ avec $j_i=0$.
	
	Comme il n'y a pas d'annulation des puissances de $x$ dans les divers termes de la forme $x^{-b_ip^i}\phi_i(x)Y$, on peut analyser chaque terme $x^{-b_ip^i}\phi_i(x)Y$ pour ses contributions au terme constant: 
	
	Cette contribution est exactement $Y_i=\begin{cases}
		 1\times \smat{1 & 0 \\ 0 & 1} \times Y & \text{si }b_i=0\\
		
		 x^{-b_ip^i}\times \smat{0 & 0 \\ 0 & x^{b_ip^i}} \times Y&   \text{si }b_i\neq0.
	
		\end{cases}$
\end{proof}

\begin{remarque}
	On peut appliquer les méthodes du théorème au cas d'un produit de Kronecker des matrices à logarithmes à plusieurs variables. Dans le cas supersingulier de deux variables, cela donne une description en termes de chiffres $p$-adiques de la matrice de logarithme de Lei \cite{leicanadian}, généralisant \cite[Section 4]{DL}.
\end{remarque}

\bibliography{distributions}{}
\bibliographystyle{plain-fr}
\end{document}

%% file: canonicalheader.tex
\usepackage{amssymb}
\usepackage{amscd,amsfonts,amsmath,amssymb, bbm,dsfont}
\usepackage{amsmath}
\usepackage{amsfonts}
\usepackage{amssymb}\usepackage{enumerate,epsf,fancyhdr,float,graphicx,tabularx}
\usepackage{latexsym,mathrsfs,multirow}
\usepackage{wasysym}
\usepackage{xypic}
\usepackage[all]{xy}\usepackage[OT2,T1]{fontenc}
\usepackage[modulo,mathlines,displaymath, running]{lineno}
\usepackage{amsmath,mathrsfs,enumerate,wasysym}
\usepackage{xypic,hhline}
\usepackage[all]{xy}
\usepackage[OT2,T1]{fontenc}
\usepackage[modulo,mathlines,displaymath]{lineno}
\usepackage{tikz,tikz-cd}
\usepackage{dashrule}
\usetikzlibrary{matrix}
\usepackage[modulo,mathlines,displaymath]{lineno}

\newtheorem{theorem}{Theorem}[section]
\DeclareSymbolFont{cyrletters}{OT2}{wncyr}{m}{n}\DeclareMathSymbol{\Sha}{\mathalpha}{cyrletters}{"58}

\renewcommand{\phi}{{\varphi}}

\newcommand{\Zp}{\mathbf{Z}_p}

\renewcommand{\geq}{\geqslant}
\renewcommand{\leq}{\leqslant}

\newcommand{\smat}[1]{\left( \begin{smallmatrix*} #1 \end{smallmatrix*} \right)}

\newcommand{\links}{\left(\begin{array}{cc}}
\newcommand{\rechts}{\end{array}\right)}
\newcommand{\bai}{\left[\begin{array}{cc}}
\newcommand{\dai}{\end{array}\right]}
\newcommand{\hidari}{\left(\begin{array}{c}}
\newcommand{\migi}{\end{array}\right)}

\newcommand{\C}{\mathbb{C}}

\newcommand{\Q}{\mathbb{Q}}

\newcommand{\Z}{\mathbb{Z}}

\newcommand{\Log}{\mathcal Log}

\renewcommand{\contentsname}{Contents\\{\footnotesize\normalfont(A table
of contents should normally not be included)}}

\newtheorem{auxiliary proposition}[theorem]{Auxiliary Proposition}

\newtheorem{main conjecture}[theorem]{Main Conjecture}
\newtheorem{main theorem}[theorem]{Main Theorem}
\newtheorem{modesty proposition}[theorem]{Modesty Proposition}

\newtheorem{open problem}[theorem]{Open Problem}

\newtheorem{déf}[theorem]{Définition}
\newtheorem{thé}[theorem]{Théorème}
\newtheorem{lemme}[theorem]{Lemme}
\newtheorem{proposition}[theorem]{Proposition}

\newtheorem{remarque}[theorem]{Remarque}

\newtheorem{convergence lemma}[theorem]{Convergence Lemma}
\newtheorem{corrected lemma}[theorem]{Corrected Lemma}
\newtheorem{growth lemma}[theorem]{Growth Lemma}
\newtheorem{coefficient lemma}[theorem]{Integrality Lemma}
\newtheorem{interpolation lemma}[theorem]{Interpolation Lemma}
\newtheorem{kernel lemma}[theorem]{Kernel Lemma}
\newtheorem{limit lemma}[theorem]{Limit Lemma}
\newtheorem{tandem lemma}[theorem]{Modesty Lemma}
\newtheorem{zero-finding lemma}[theorem]{Zero-Finding Lemma}